\providecommand{\href}[2]{#2}

\documentclass[12pt]{amsart}

\usepackage{graphicx}
\usepackage{amsthm}
\usepackage{amsmath, amssymb}
\swapnumbers

\theoremstyle{plain}
\newtheorem{Thm}{Theorem}

\newtheorem{Coro}[Thm]{Corollary}
\newtheorem{Lem}[Thm]{Lemma}

\theoremstyle{definition}
\newtheorem{Def}[Thm]{Definition}

\begin{document}

\title{ Trisections and link surgeries}

\maketitle

 \author{Robion Kirby and Abigail Thompson}  \footnote{Both authors are supported in part by the National Science Foundation.}

\begin{abstract} {We examine questions about surgery on links which arise naturally from the trisection decomposition of 4-manifolds developed by Gay and Kirby \cite{G-K3}.  These links lie on Heegaard surfaces in $\#^j S^1 \times S^2$ and have surgeries yielding $\#^k S^1 \times S^2$.   We describe families of links which have such surgeries.     One can ask whether all links with such surgeries lie in these families;  the answer is almost certainly no.   We nevertheless give a small piece of evidence in favor of a positive answer}.     
 \end{abstract}

\section{Introduction} 

The question of which manifolds can arise from Dehn surgery on a knot in the 3-sphere is much-studied,   from Gabai's proof of Property R \cite{G}, to Gordon and Luecke's solution to the knot complement problem \cite{GL}, to the as-yet open question posed by the Berge conjecture \cite{K} regarding which knots have surgeries yielding lens spaces.   Expanding our attention to consider multi-component links in the 3-sphere, even to 2-component links,  has proven remarkably difficult.  For example, attempting the most straightforward generalization of Property R to a 2-component link in the 3-sphere (i.e., characterizing 2-component links with surgeries yielding  $\#^2 S^1 \times S^2$) has given rise to potential counter-examples to the slice-ribbon conjecture \cite{GST}.    We suggest a framework inspired by trisections of 4-manifolds into which many of these questions can be placed.         

In the next section we review basic terminology of trisections and Heegaard splittings.    
In section 3 we introduce some new definitions and use them to prove two lemmas about link surgeries in arbitrary 3-manifolds yielding some number of copies of $\#S^1 \times S^2$.   In section 4  we pose some questions in the specific context arising from trisections, and in section 5 we prove our main result.   

\section{Background and definitions} 

Let $X$ be a closed, orientable, smooth 4-manifold.    In \cite{G-K3} Gay and Kirby show that $X$ has a trisection into three 4-dimensional handlebodies, and prove that any two trisections of $X$ are stably equivalent under a suitable notion of stabilization.   

\begin{Def} \label{D:trisection}
 A \emph{$(g;k_1,k_2,k_3)$--trisection} of a closed, oriented $4$--manifold $X$ (where $0 \leq k_i \leq g, i=1,2,3$) is a decomposition $X = X_1 \cup X_2 \cup X_3$, where (1) each $X_i \cong \natural^{k_i} S^1 \times B^3$, (2) each $X_i \cap X_j \cong \natural^g S^1 \times B^2$ (for $i \neq j$), and (3) $X_1 \cap X_2 \cap X_3 \cong \#^g S^1 \times S^1$.
\end{Def}

The topology of each of the three pieces of $X$ is completely determined by a single integer $k_i$, and the topology of each of the overlaps between pieces is determined by another integer $g$.  If  $k=k_1=k_2=k_3$ the trisection is called {\em balanced}.   We are particularly interested in balanced trisections with $k=0$.

Given a trisection of $X^4$, we have a {\em central surface} $\Sigma=X_0 \cap X_1 \cap X_2$ in $X$ bounding three $3$--dimensional handlebodies $X_i \cap X_j$, which fit together in pairs to form Heegaard splittings of three $3$--manifolds in $X$, and these $3$--manifolds in turn uniquely bound three $4$--dimensional $1$--handlebodies.   We can thus specify a trisection by considering systems of curves on $\Sigma$:

\begin{Def}
 A {\em cut system} for a closed surface $\Sigma$ of genus $g$ is a collection of $g$ disjoint simple closed curves on $\Sigma$ which cut $\Sigma$ open into a $2g$--punctured sphere.
\end{Def}

\begin{Def}
 A genus $g$ {\em Heegaard diagram} for a closed orientable 3-manifold is a triple $(\Sigma, \alpha, \beta)$, where $\Sigma$ is a closed orientable genus $g$ surface and each of $\alpha$ and $\beta$ is a cut system for $\Sigma$. 
 
 \end{Def}

\begin{Def}
A {\em $(g;k_1,k_2,k_3)$--trisection diagram} is a $4$--tuple $(\Sigma, \alpha, \beta, \gamma)$ such that each of $(\Sigma,\alpha,\beta)$, $(\Sigma, \beta, \gamma)$ and $(\Sigma,\gamma,\alpha)$ are genus $g$ Heegaard diagrams of $\#^k{_i} S^1 \times S^2, i=1,2,3$ respectively. A trisection diagram for a given trisection $X=X_1 \cap X_2 \cap X_3$ is a trisection diagram $(\Sigma, \alpha, \beta, \gamma)$, where $\Sigma$ is diffeomorphic to $X_1 \cap X_2 \cap X_3$, $\alpha$ is a cut system for $X_1 \cap X_2$, $\beta$ for $X_2 \cap X_3$, and $\gamma$ for $X_3 \cap X_1$.
\end{Def}

\section{Link surgeries in $N^3$ yielding $\#^k S^1 \times S^2$}

In the case of a balanced trisection with $k=0$, each of the pairs of cut systems define a copy of $S^3$, and we can think of building $M^4$ in the following way:

\noindent Start with a 4-ball with boundary $S^3$.\\
Do framed surgery on the $\gamma$ curves to obtain $S^3$ again.\\
Cap off the resulting object with another 4-ball.

The collection of $g$ $\gamma$ curves in $S^3$ are a link $L$ lying on a genus g Heegaard surface $\Sigma$ for $S^3$. $L$ is a cut system on $\Sigma$.

\begin{figure}[h]
    \centering
    \includegraphics[width=0.4\textwidth]{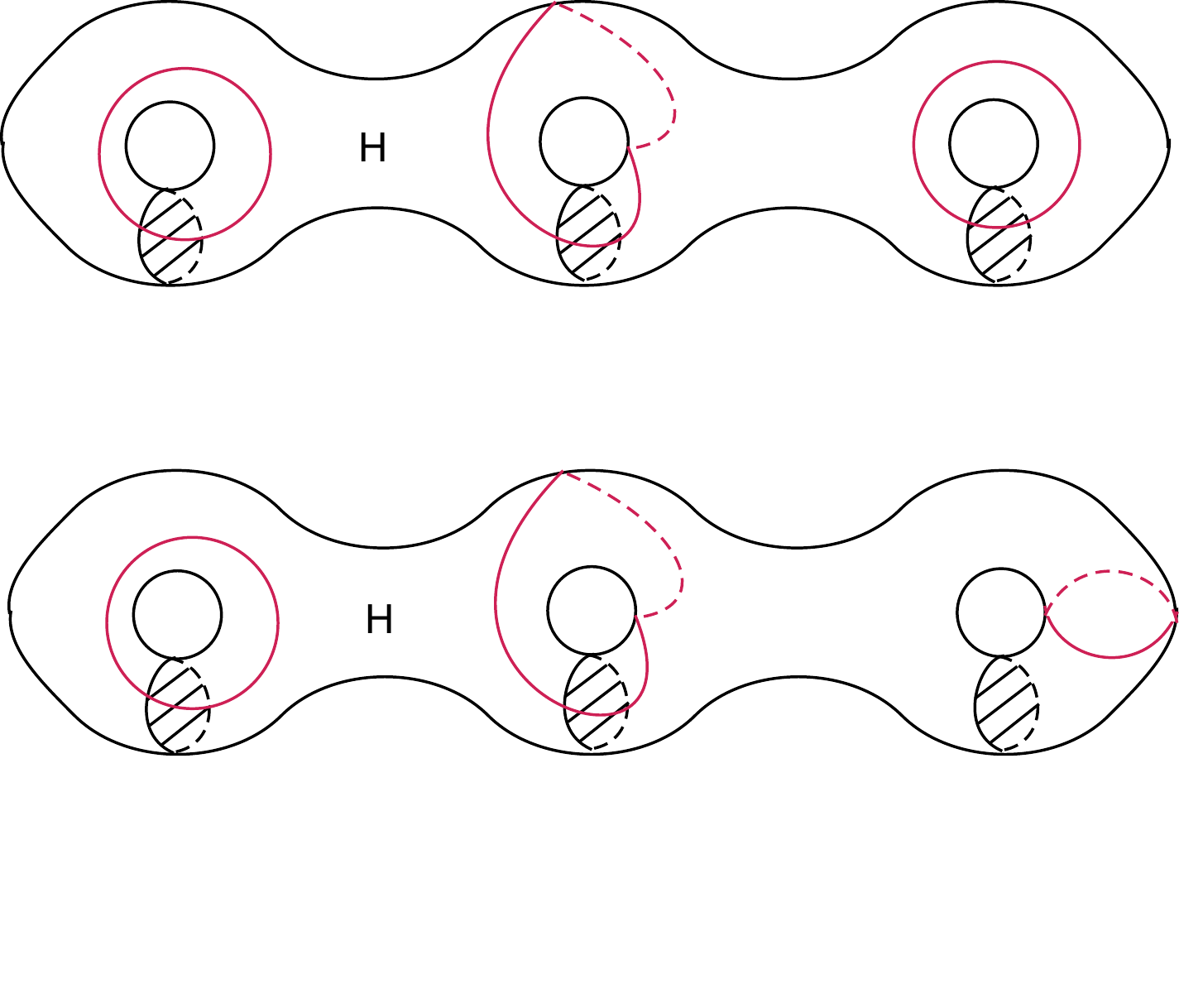}
    \caption{primitive system (top) and pseudo-primitve system (bottom)}
    \label{primitive}
\end{figure}

What follows are some definitions and observations using this set-up as an inspiration:

Let $L$ be a g-component link imbedded as a cut system on the boundary of a genus g handlebody $H$.   We say that $L$ is {\em primitive} on $H$ if there exists a complete set of compressing disks for $H$ whose boundaries are geometrically dual to $L$ (see figure \ref{primitive}, top).   Note that the boundary of $H$ naturally induces a framing on the components of $L$, the {\em surface framing}.  We say that $L$ is {\em slide-primitive}  on $H$ if it is possible to do surface-framed handle slides on $L$ lying within $\partial{H}$ to obtain a link $L'$ such that $L'$ is primitive on $H$.

Let $L$ be a g-component link imbedded as a cut system on a genus g Heegard surface $\Sigma$ in a 3-manifold $N^3$, with $\Sigma$ bounding handlebodies $H_1$ and $H_2$.    

We say that $L$ is {\em double-slide-primitive} (or{ \em dsp}) on $\Sigma$ if $L$ is slide-primitive on $H_1$ and on $H_2$.

\begin{Lem}\label{Lemma1}
Let $L$ be a g-component link imbedded as a cut system on a genus g Heegard surface $\Sigma$ in a 3-manifold $N^3$, with $\Sigma$ bounding handlebodies $H_1$ and $H_2$.     Then surface-framed surgery on $L$ yields $S^3$ iff $L$ is dsp on $\Sigma$.
\end{Lem}

Proof:

Assume surgery on $L$ yields a 3-sphere.   The cut system $L$ caps off $\Sigma-L$ into a 2-sphere after surgery on $L$.    Since the result of surgery is $S^3$, this 2-sphere must bound a 3-ball on both sides.     Thus each of the manifolds $M_i=H_i\cup$(2-handles), where the 2-handles are attached along $L$, must be 3-balls.   The result follows from Waldhausen's theorem \cite{Wal} on Heegaard splittings of the 3-sphere applied to each $M_i$.

Conversely if $L$ is dsp on $\Sigma$, then each of the $M_i$ are 3-balls and so the surgered manifold is the 3-sphere.

\bigskip

We can generalize this:\\

Let $L$ be a g-component link imbedded as a cut system on the boundary of a genus g handlebody $H$.   We say that $L$ is {\em pseudo-primitive} on $H$ if there exists a complete set of compressing disks for $H$ whose boundaries are geometrically dual to $L$ or isotopic to curves in $L$.   We say that $L$ is {\em slide-pseudo-primitive}  on $H$ if it is possible to do surface-framed handle slides on $L$ to obtain a link $L'$ such that $L'$ is pseudo-primitive on $H$ (see figure \ref{primitive}, bottom).

Let $L$ be a g-component link imbedded as a cut system on a genus g Heegard surface $\Sigma$ in a 3-manifold $N^3$, with $\Sigma$ bounding handlebodies $H_1$ and $H_2$.     

We say that $L$ is {\em double-slide-pseudo-primitive} (or{ \em dspp}) on $\Sigma$ if $L$ is slide-pseudo-primitive on $H_1$ and on $H_2$.

\begin{Lem}\label{Lemma2}
Let $L$ be a g-component link imbedded as a cut system on a genus g Heegard surface $\Sigma$ in a 3-manifold $N^3$, with $\Sigma$ bounding handlebodies $H_1$ and $H_2$.    Then surface-framed surgery on $L$ yields $\#^k S^1 \times S^2$ iff $L$ is dspp on $\Sigma$.
\end{Lem}

Proof:

Assume surgery on $L$ yields $\#^k S^1 \times S^2$.   The cut system $L$ caps off $\Sigma-L$ into a 2-sphere after surgery on $L$.    Since the result of surgery is $\#^k S^1 \times S^2$, this 2-sphere must bound $k_i$, $i=1,2$, copies of $\# S^1 \times S^2$  on both sides, where $k_1+k_2=k$.     Thus each of the manifolds $M_i=H_i\cup$(2-handles), where the 2-handles are attached along $L$, is $\#^{k_i} S^1 \times S^2$.   The result follows from the generalization of Waldhausen's theorem ( \cite{SS}, p. 313) on Heegaard splittings of the 3-sphere to Heegaard splittings of $\#^{k_i} S^1 \times S^2$.

Conversely if $L$ is dspp on $\Sigma$, then each of the $M_i$ are $\#^{k_i} S^1 \times S^2$ and so the surgered manifold is $\#^k S^1 \times S^2$.    \\

In the context of trisections of 4-manifolds, $N^3$ in Lemma \ref{Lemma2} is always itself some number of copies, say $j$,  of $S^1 \times S^2$.   In the next section we consider this case specifically.   

\section{Link surgery questions arising from trisections}

Suppose $L$ is a g-component link imbedded as a cut system on a genus g Heegard surface in $\#^j S^1 \times S^2$ and surface-framed surgery on $L$ yields $\#^k S^1 \times S^2$.   Allowing arbitrary handle slides on $L$ will not change the result of the surgery, but may turn $L$ into a link $L'$ which no longer lies on the Heegaard surface.   One could conjecture that this is essentially the only way to generate such an $L$.           

Even more optimistically, in the special case where $N^3$ is the 3-sphere (so $j=k=0$)  and surgery on $L$ yields the 3-sphere back, one could hope $L$ is even simpler than slide-equivalent to a ``dsp" link on a Heegaard surface, and ask the following:\\

\noindent {\em Question 1}:   Let $L$ be a g component framed link in the 3-sphere such that surgery on $L$ yields the 3-sphere.     Is $L$ handle-slide equivalent to a union of Hopf links and unknots?  \\   

The answer to this is ``no"; Harer, Kas and Kirby (\cite{HKK}, p. 66) give a handle decomposition of  $K3$  which utilizes a 22-component link $L$ in the 3-sphere with integral surgery yielding $S^3$ which cannot be handle-slide equivalent  to a union of Hopf links and unknots ($K3$ has signature 16 but Hopf links and unknots give signature 0).    Meier and Lambert-Cole \cite{MLC} have noted that this decomposition corresponds to a genus 22 balanced trisection of $K3$ with all sectors a 4-ball.    In particular the Akbulut-Kirby link lies on a genus 22 Heegaard surface in the 3-sphere.    $L$ is visibly dsp on this Heegaard surface (as Lemma \ref{Lemma1} requires).       

We increase our chance of success by asking for less:\\

\noindent {\em Question 2}: Let $L$ be a g component framed link in the 3-sphere, such that surgery on $L$ yields the 3-sphere.     Must $L$ be handle-slide equivalent to a dsp link on a genus g Heegaard surface for $S^3$?  \\

The answer to this is likely also ``no", but we don't know any counter-examples.   A small piece of evidence in the ``yes'' direction for this question is our main theorem, which will appear in the final section. 

To conclude this section, we state the most general form of the surgery question as it arises in the trisection context, noting a result which provides support for a negative answer:\\

\noindent {\em Question 3}: Let $L$ be a $g\geq{k}$ component framed link in  $\#^k S^1 \times S^2$ such that surgery on $L$ yields $\#^j S^1 \times S^2$.     Must $L$ be handle-slide equivalent to a dspp link on a genus g Heegaard surface for $\#^k S^1 \times S^2$? \\

Again presumably the answer is ``no"; an example of why this presumption is justified can be found in \cite{GST}.   There, a 2-component link $L$ in $S^3$ (actually a whole family of 2-component links) is shown to yield $\#^2 S^1 \times S^2$ after 0-framed surgery.   If $L$ were handle-slide equivalent to a dspp link on a genus 2 Heegaard surface for $S^3$, it would be handle-slide equivalent to an unlink (\cite{GST}, proposition 3.1).    This would imply that conjectured counter-examples to the Andrews-Curtis conjecture fail.

\section{Main theorem}

\begin{Thm}\label{Main}
Let $L=L_1\cup{L_2}$ be a framed 2-component link in the 3-sphere such that surgery on $L$ yields $S^3$.   Suppose $L_1$ is the unknot.\\
Then $L$ is handle-slide equivalent to a dsp link on a genus 2 Heegaard surface for $S^3$.   
\end{Thm}

A stronger statement follows immediately using a result of Meier and Zupan \cite{M-Z}:\\
\begin{Coro}\label{G2}
Let $L=L_1\cup{L_2}$ be a framed 2-component link in the 3-sphere such that surgery on $L$ yields $S^3$.   Suppose $L_1$ is the unknot.\\
Then $L$ is handle-slide equivalent to the unlink or Hopf link.   

\end{Coro} 

Proof of Corollary \ref{G2}:

Meier and Zupan \cite{M-Z} classify all genus 2 trisection diagrams.  The corollary follows from this classification once $L$ has been handle-slid using Theorem \ref{Main}  to lie on a genus 2 Heegaard surface for $S^3$.\\

\begin{figure}[hbtp]
    \centering
    \includegraphics[width=0.3\textwidth]{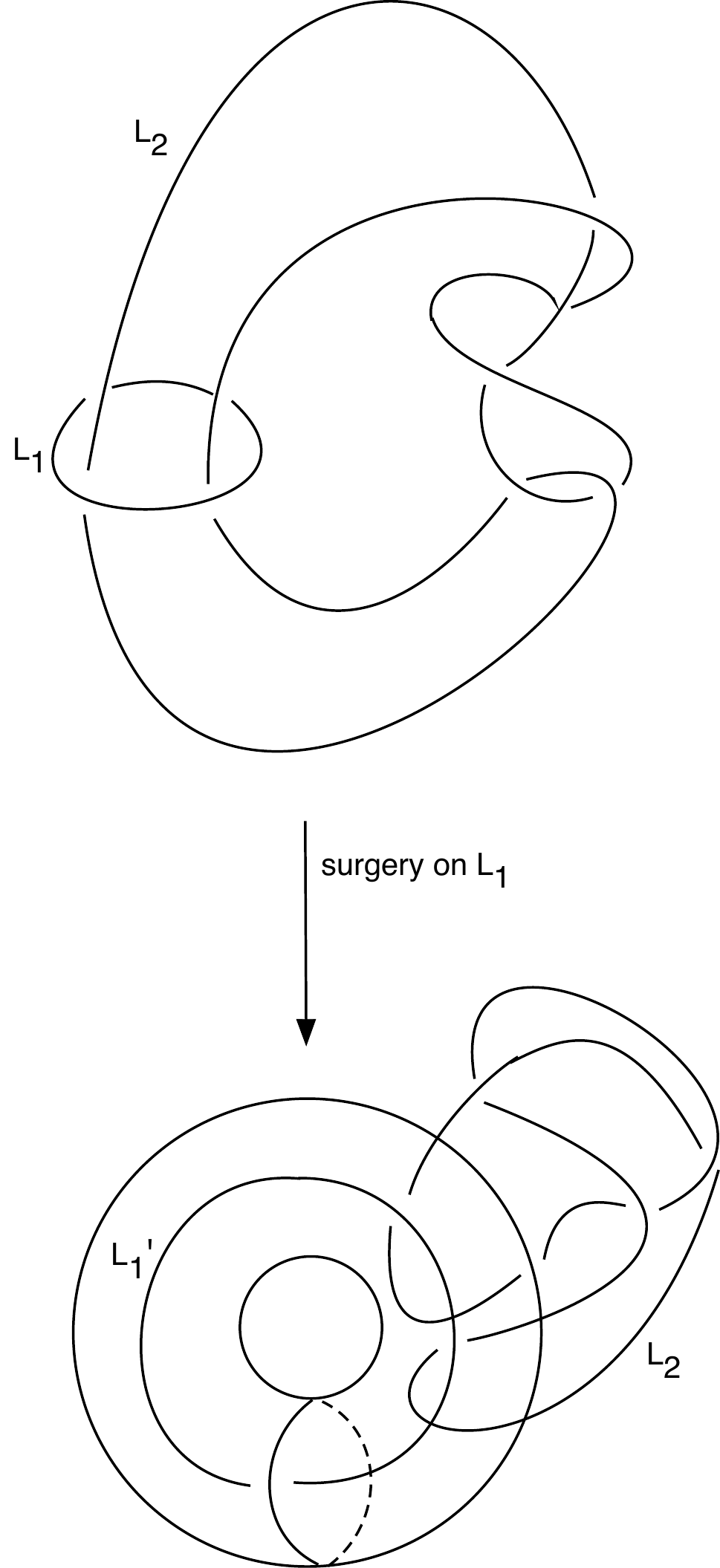}
    \caption{}
    \label{surg2}
\end{figure}

\begin{figure}[hbtp]
    \centering
    \includegraphics[width=0.3\textwidth]{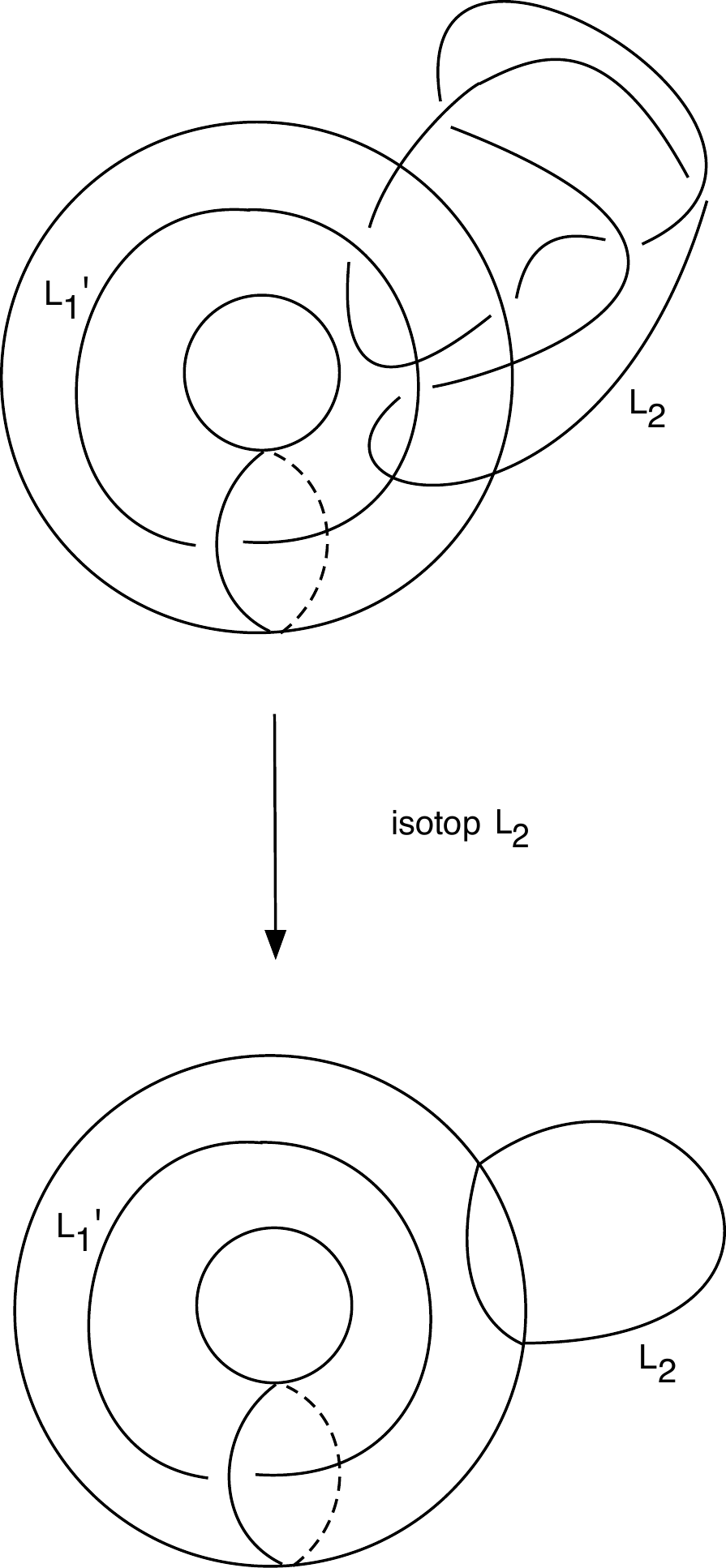}
    \caption{}
    \label{surg3}
\end{figure}

Proof of Theorem \ref{Main}:

Do the indicated surgery on $L_1$ (see Figure \ref{surg2}).   Since $L_1$ is the unknot, this yields $M=S^3$, $S^1\times{S^2}$ or a lens space L(n,1).    $L_2$ is then a knot in $M$, and surgery on $L_2$ must yield $S^3$.   If $M$ is $S^3$ or $S^1\times{S^2}$, $L_2$ is respectively the unknot \cite{GL} or a core curve of a genus one Heegaard splitting \cite{G}.  The trickiest case is the lens space case.   However we are considering only integer surgeries, and the Berge conjecture \cite{K} is known in L(n,1) \cite{KMOS}.  Hence in all cases  $L_2$ can be isotoped in $M$ to have bridge number zero or  one with respect to the genus one Heegaard splitting defined by $L_1$ (see Figure \ref{surg3}).     During this isotopy, $L_2$ may cross the dual knot $L_1'$.   Each such crossing corresponds to a handle-slide of $L_2$ over $L_1$ in the original link diagram (see the proof of proposition 3.2 in  \cite{GST}).     Once the isotopy is complete (or the corresponding handle slides in the 3-sphere are complete), the new link $L_1\cup{L_2'}$ is obviously tunnel number one, and so the link $L'=L_1\cup{L_2'}$ in $S^3$ can be imbedded as required (i.e., as a cut system with the desired framing) on a genus 2 Heegaard surface in $S^3$.  The theorem  follows from Lemma \ref{Lemma1}.          \\

\

\section{Acknowledgement}  We thank Jeff Meier for many helpful comments resulting in considerable simplification of the main result.

\end{document}